\documentclass{amsart}
\usepackage{amsfonts,amssymb,amscd,amsmath,enumerate,verbatim,calc}

\newcommand{\CM}{Cohen-Macaulay}

\newcommand{\ff}{\text{if and only if}}
\newcommand{\wrt}{with respect to}
\newcommand{\B}{\mathfrak{b} }

\newcommand{\m}{\mathfrak{m} }
\newcommand{\Pf}{\mathfrak{P} }

\newcommand{\A}{\mathfrak{a} }

\newcommand{\Z}{\mathbb{Z} }

\newcommand{\xar}{\longrightarrow}

\newcommand{\reg}{\operatorname{reg}}

\newcommand{\gr}{\operatorname{gr}}

\newcommand{\height}{\operatorname{height}}

\newcommand{\projdim}{\operatorname{projdim}}
\newcommand{\Spec}{\operatorname{Spec}}
\newcommand{\Min}{\operatorname{Min}}
\newcommand{\Ass}{\operatorname{Ass}}
\newcommand{\Assh}{\operatorname{Assh}}

\newtheorem{theorem}{Theorem}[section]

\newtheorem{conjecture}[theorem]{Conjecture}
\theoremstyle{definition}

\theoremstyle{remark}
\newtheorem{remark}[theorem]{Remark}
\newtheorem{s}[theorem]{}
\newtheorem{observation}[theorem]{Observation}
\numberwithin{equation}{section}

\begin{document}

\title[Multiplicity Conjecture]{On the upper bound of \\ the Multiplicity conjecture}
\author{Tony~J.~Puthenpurakal}

\address{Department of Mathematics, IIT Bombay, Powai, Mumbai 400 076}

\email{tputhen@math.iitb.ac.in}

\thanks{ The author thanks Universt{\"a}t
Duisburg-Essen 
  for hospitality during Nov-Dec 2006. The author thanks DFG for financial support, which made this visit possible}

\subjclass[2000]{Primary  13H15, 13D02 ; Secondary 13D40, 13A30}


\dedicatory{Dedicated to Juergen Herzog on the occasion of his
65th birthday}

\keywords{multiplicity conjecture, regularity,  reduction, analyticity}

\begin{abstract}
Let $A = K[X_1,\ldots,X_n]$ and let $I$ be a graded ideal in $A$. We show that the upper bound of Multiplicity conjecture of Herzog, Huneke and Srinivasan holds asymptotically (i.e., for $I^k$ and all $k \gg 0$) if $I$ belongs to any of  the following large classes of ideals:
\begin{enumerate}[\rm (1)]
 \item 
 radical ideals.
\item 
 monomial ideals with generators in different degrees.
\item
 zero-dimensional ideals with generators in different degrees.
\end{enumerate}
Surprisingly, our proof uses local techniques like analyticity, reductions, equimultiplicity and local results like Rees's theorem
on multiplicities.
\end{abstract}

\maketitle
\section{introduction}
Let $K$ be a field and  $A = K[X_1,\ldots,X_n]$ be a polynomial ring with standard grading. Let $I$ be a graded ideal of $A$.  Let
$$
0
\xar
\bigoplus_{j\in \Z} A(-j)^{\beta_{p,j}(A/I)}
\xar
\cdots
\xar
\bigoplus_{j\in \Z} A(-j)^{\beta_{1,j}(A/I)}
\xar A
\xar 0
$$
be a minimal graded free resolution of $A/I$. Set $ p = \projdim A/I$ and $c = \height I$.
Consider for $1\leq i \leq p$ the numbers
$$
M_i( A/I) = \max\{j \in \Z \mid \beta_{i,j}(A/I)\neq 0\}
\ \text{\&}  \ 
m_i(A/I) = \min\{j \in \Z \mid \beta_{i,j}(A/I)\neq 0\}.
$$
Let $e(A/I)$ denote the multiplicity of $A/I$. Set
$$L(I) =  \frac{1}{c!}\prod_{i=1}^{c} m_i(A/I)  \quad \text{ and } \quad U(I) = \frac{1}{c!}\prod_{i=1}^{c} M_i(A/I).$$

The  conjecture of Herzog, Huneke and Srinivasan states that
\begin{conjecture}\label{C-cohen}
If $A/I$ is \CM \ then
\[ 
L(I) \leq e(A/I) \leq U(I).
 \]
\end{conjecture}
If $A/I$ is not \CM \ then in \cite{HS} it is conjectured that
\begin{conjecture}\label{C-gen}
\[ 
 e(A/I) \leq U(I).
 \]
\end{conjecture}
Both Conjectures \ref{C-cohen} and \ref{C-gen} have been proved for many classes of ideals  (see \cite{GO03, GUVT05,  HS, HEZH, MINARO05, MIR06, NOSC,RO05,SR98}). For extensions of this conjecture see \cite{ MINARO06, MINAZA06, MIR07,RO07, ZA05}. For some new approaches to this problem see
 \cite{BOSO, FR06,  GOSCSR05, MINARO06, MINAZA06}. Our result is

\begin{theorem}\label{main}
 Let $I$ be a graded ideal. If $I$ belongs to any of the following classes of ideals
\begin{enumerate}[\rm (1)]
 \item 
radical ideals.
\item 
 monomial ideals with generators in different degrees.
\item
 zero-dimensional ideals with generators in different degrees.
\end{enumerate}
Then $e(A/I^k) \leq U(I^k)$ for all $ k \gg 0.$
\end{theorem}
In \cite[Theorem 2]{HEZH},  the authors show $\lim_{k\to
\infty} e(A/I^k) /U(I^k) \leq 1$. There are examples where  $\lim_{k\to
\infty} e(A/I^k) /U(I^k) = 1$, for instance see Section 4.
In our proof we show that  in the class of ideals of Theorem \ref{main} we have that the limit on the left hand side
is $< 1$. The surprising feature of our proof is the use of local techniques like equimultiplicity, reductions analyticity and local theorems like Rees multiplicity theorem (see \ref{rees}).

\textit{Overview of  the paper.} In section two we introduce notation and discuss a few preliminary facts that we need. In section three we prove Theorem \ref{main}. In section 4 we give an example of a class of ideals which satisfy $\lim_{k\to
\infty} e(A/I^k) /U(I^k) = 1$.

\textit{Acknowledgment:} Its a pleasure to thank Prof. J. Herzog  for many discussions regarding this paper. I also thank the referee for careful reading. 

\section{Preliminaries}
In this section we recall some notions in local algebra. We also discuss asymptotic behavior of regularity of ideals
$I^k$ for $k \gg 0$. Finally we also recall that the function $k \mapsto e(A/I^k)$ is polynomial in $k$ for $k \gg 0$.

$\bullet$ \textit{Some local notions:}

Let $(R,\m)$ be a Noetherian local ring of dimension $d$ and residue field $K = R/\m$ which, for
convenience, we assume is infinite.   Let $\A$ be an ideal in $R$. If $M$ is a finitely generated $R$-module then
$\mu(M)$ denotes its  minimal number of generators and $\ell(M)$ denotes its length. 

\s\label{spread} The \textit{analytic spread}
of $\A$ is the Krull dimension of the fiber-cone $F(\A) = \bigoplus_{n \geq 0}\A^n/\m \A^n$. We denote it by $s(\A)$.
By \cite[p.\ 150, Th. 1]{NoR}, $s(\A) = \mu(\B)$ where $\B$ is a (any) minimal reduction  of $\A$. For definition of  reduction
and minimal reduction see  \cite[p.\ 146]{NoR}.
It can be shown that $\height(\A) \leq s(\A)$, see \cite[p.\ 151, L. 4]{NoR}. We say $\A$ is an \textit{equimultiple} ideal if $\height(\A) = s(\A)$. 
If $R$ is quasi-unmixed then $\A$ is equimultiple \ff \  $\gr_{\A} R = \bigoplus_{n \geq 0}\A^n/\A^{n+1}$, the associated graded ring of $\A$, has a homogeneous system of parameters, see \cite[2.6]{GHO}.

\s\label{rees} If $\A$ is $\m$-primary then let $e(\A, R) = $  multiplicity of $R$ \wrt \ $\A$ i.e.,
\[
 e(\A, R) = \lim_{n \rightarrow \infty} \frac{d!}{n^d}\ell\left(\frac{R}{\A^n} \right).
\]
Let $\B \subseteq \A$ be $\m$-primary. Clearly $e(\B, R) \geq e(\A, R)$. It is easy to see that if $\B$ is a reduction of $\A$ then $e(\B, R) = e(\A, R)$. A celebrated theorem due to Rees \cite{Rees1}
shows that  if $R$ is quasi-unmixed  and  $e(\B, R) = e(\A, R)$  then $\B$ is a reduction of $\A$.

$\bullet $\textit{ Asymptotic behavior of regularity:}

Let $A = K[X_1,\ldots,X_n]$. Let $I$ be a  graded ideal in $A$ and let $p = \projdim A/I$.
Then $$\reg(I)=\max\{ M_{i+1}(A/I) - i \mid i=0,\ldots, p-1\}$$ is the regularity of
$I$.
Set $\reg_i(I) = M_{i+1}(A/I) - i$ for  $i=0,\ldots, p-1$.

\s\label{basicreg} 
In \cite[2.4]{CHT} and \cite[1]{K}  it is shown that $\reg(I^k) = qk + r$ for $k \gg 0$. In \cite[5]{K} it is shown that $I$ has a reduction $J$ such that $\reg_0(J) = q$.  In particular
 $J$ is generated in degrees $\leq q$.  We call such a reduction to be a \textit{Kodiyalam reduction.}

\s \label{zheng}
  In \cite[2.1(ii)]{HEZH}
it is proved that for $i = 0,\ldots, c-1$, 
\[
\reg_i(I^k) = qk + r_i \quad \text{ for all} \ \ k \gg 0.
\]

Therefore \text{for} $ k\gg 0$,
\[
 U(I^k) = \frac{1}{c!}\prod_{i=1}^{c} M_i(A/I^k)=\frac{q^c}{c!}k^c+\cdots\quad +\text{lower terms in $k$} 
\]

$\bullet$.\textit{ The function} $k \mapsto e(A/I^k)$.
\s \label{function} Let $\Assh(I) = \{ \Pf\in \Spec(A) \mid  \Pf \supseteq I \ \text{and} \ \dim A/\Pf = \dim A/I \}$. Notice that all $\Pf \in \Assh(I)$ are 
graded ideals.  The associativity
formula of multiplicity (\cite[4.7.8]{BH}) then shows
that
\begin{equation}\label{z1}
 e(A/I)=\sum_{\Pf \in \Assh(I)}\ell(A_{\Pf}/I_{\Pf})e(A/\Pf).
\end{equation}
Since  $\Assh(I^k) = \Assh(I)$ for all $k \geq 1$ we have that 
\begin{equation}\label{formMult}
 e(A/I^k)=\sum_{\Pf \in \Assh(I)}\ell(A_{\Pf}/I_{\Pf}^k)e(A/\Pf).
\end{equation}
Recall that $c = \height(I)$. Since $k \mapsto \ell(A_{\Pf}/I_{\Pf}^k)$ is a polynomial function of degree $c$ it follows that $ k \mapsto e(A/I^k)$ is a polynomial
function of degree $c$. Furthermore if $E(I)$ is the normalized leading coefficient of this function then
\begin{equation}\label{z3}
 E(I) = \sum_{\Pf \in \Assh(I)}e(I_\Pf , A_{\Pf})e(A/\Pf).
\end{equation}

\begin{remark}\label{e-red}
Let $J \subseteq I$ be a graded ideal. If $J$ is a reduction of $I$ then
$J_{\Pf}$ is a reduction of $I_{\Pf}$ for all primes $\Pf$. So
  $E(J) = E(I)$. 
\end{remark}

\s \label{formula} By \ref{zheng} and \ref{function} we get that 
\[
 \lim_{k \mapsto \infty} \frac{e(A/I^k)}{U(I^k)}  = \frac{E(I)}{q^c}
\]
Here $q$ is as in \ref{basicreg}.
\section{Proof of Theorem \ref{main}}
In this section we prove our result. We use 
 \cite[Theorem 2]{HEZH}, where it is proved that $\lim_{k\to
\infty}e(A/I^k)/U(I^k)  \leq 1$.
In our proof we show that  in the class of ideals of Theorem \ref{main} we have that the limit on the left hand side
is $< 1$. Throughout this section $q$ is as in \ref{basicreg}.

 \s \label{H-Z}
In  \cite[section 2]{HEZH} the authors  assume $K$ is infinite and then do the following:
\begin{itemize}
 \item 
Let  $J$ be a Kodiyalam reduction of $I$ and let $f_1, \ldots, f_c \in J_q$ be $c$-generic $q$-forms. Set
$L = (f_1,\ldots, f_c)$.
\item 
$ e(A/L) = E(L) = q^c.  $
\item 
$E(I) \leq E(L).$
\end{itemize}
The folllowing observation is useful:
\begin{observation}\label{obs}
\begin{enumerate}[\rm (1)]
 \item 
In \ref{H-Z}, an ideal $L =(f_1,\ldots,f_c)$ where $f_1, \ldots, f_c \in J_q$ is a regular sequence will do. In fact 
in \cite[section 2]{HEZH} it is chosen generic just to ensure $f_1, \ldots, f_c $ is a regular sequence.
\item 
We may choose $f_1 \in J_q$ to be any non-zero element.
\end{enumerate}
\end{observation}
To prove $E(I) < E(L)$
the following remark is useful:
\begin{remark}\label{ineq}
 $L$ is unmixed. Also $\height I = \height L = c$. Thus $\Assh(I) \subseteq \Assh(A/L) = \Min(A/L)$; the set of minimum primes of $L$. So to prove $E(I) < E(L)$, it suffices to show that there exists $\Pf \in \Assh(I)$ such that 
$e(I_{\Pf},A_{\Pf}) < e(L_{\Pf}, A_{\Pf})$.
\end{remark}
We now give
\begin{proof}[Proof of Theorem \ref{main}]
We prove that for each of the class of ideals considered we have $E(I) < E(L) = q^c$. We also assume $K$ is infinite.
This follows from the usual standard trick in the case when $K$ is finite. 

\textbf{Case 1:}  \textit{$I$ is a radical ideal.}
In this case we  \newline
\textit{Claim:}  $E(I) = e(A/I).$ \newline
Let $I = Q_1\cap \cdots\cap Q_s$ be a minimal irredundant primary decomposition of $I$. Set $\Pf_i = \sqrt{Q_i}$ for $i = 1,\ldots, s$.. Let $\Pf \in \Assh(I)$.
Then $\Pf = \Pf_i$ for some $i$.  
\textit{As $I$ is a radical ideal} we have $I_{\Pf} = \Pf A_{\Pf}$. Notice $A_{\Pf}$ is a regular local ring of dimension
$c$. So
\begin{alignat*}{2}
 e(I_{\Pf}, A_{\Pf}) &= e(\Pf A_{\Pf}, A_{\Pf}) &= 1, \\
\ell(A_{\Pf}/ I_{\Pf}) &= \ell(A_{\Pf}/\Pf A_{\Pf}) &=1.
\end{alignat*}
Therefore by \ref{function}(1) and (3) we get
\begin{alignat*}{2}
 E(I) &= \sum_{\Pf \in \Assh(I)}e(I_\Pf , A_{\Pf})e(A/\Pf) & \qquad  &= \sum_{\Pf \in \Assh(I)}e(A/\Pf)  \\
e(A/I) &= \sum_{\Pf \in \Assh(I)}\ell(A_{\Pf}/I_{\Pf})e(A/\Pf)   & \qquad  &= \sum_{\Pf \in \Assh(I)}e(A/\Pf)
\end{alignat*}
Thus $ E(I) = e(A/I)$. 

Set $\m = (X_1,\ldots, X_n)$. Notice that $I_\m$ is a radical ideal of $A_\m$.

\textbf{Subcase $1:$} \textit{The ideal $I_\m$ is  equimultiple}. Then by a result due to Cowsik and Nori \cite{CN}
we have that $I_\m$ is generated by a regular sequence. Since $I$ is graded it follows that 
$I$ is also generated by a regular sequence. 
In this case by \cite{GUVT05} we have that
$e(A/I^k) \leq U(I^k)$ for all $k \geq 1$. 

\textbf{Subcase $2:$}\textit{ $I_\m$ is not equimultiple.} Let $L$ be as in \ref{H-Z}. Then $L$ is not a reduction of $I$. (Otherewise $L_\m$ will be a reduction of $I_\m$ and this will imply that $I_\m$
is equimultiple.)

In particular $L \neq I$. Consider the exact sequence:
\[
 0 \xar \frac{I}{L} \xar \frac{A}{L} \xar \frac{A}{×I} \xar 0.
\]
Since $I \neq L$ we have that $ (I/L)_{\Pf} \neq 0$ for some $\Pf \in \Ass(A/L) = \Min(A/L)$. In particular
$\dim I/L =  \dim A/L = \dim A/I$. It follows that
\[
 E(L) = e(A/L) = e(A/I) + e(I/L) >  e(A/I) = E(I).
\]
This implies the result in this case. 

\textbf{Case 2:}  \textit{$I$ is a monomial ideal with generators in different degrees.} \newline
Let $\Pf \in \Assh(I)$. As $I$ is a monomial ideal, $\Pf$ is generated by a subset of the variables \cite[4.4.15]{BH}.
Say $\Pf = (X_{i_1}\ldots X_{i_s})$. Let $G(I) = \{ u_1, \ldots,u_a \}$ be the unique set of minimal
monomial generators of $I$. Assume $\deg u_1 < q$. Set $\alpha = q - \deg u_1$. 
Set $f_1 = X_{i_1}^{\alpha}u_1$ and let $L = (f_1,f_2, \ldots, f_c)$ (see \ref{obs}(2)). As $f_1 \in \Pf I_{\Pf}$,
it follows that $L_{\Pf}$ is not a minimal reduction of $I_{\Pf}$ \cite[Lemma 2]{NoR}.  Therefore by Rees's theorem
$e(L_{\Pf}, A_{\Pf}) > e(I_{\Pf}, A_{\Pf})$. So by \ref{ineq} we get $E(L) > E(I)$.
 
\textbf{Case 3:}  \textit{$I$ is a zero-dimensional ideal with generators in different degrees.} \newline
Notice that in this case $ \Assh(I) = \{(X_1,\ldots,X_n) \}$. The proof is similar to case 2.
\end{proof}

\begin{remark}
 For cases 2 and 3 in our theorem note that the ideal can never have a pure resolution. Notice also that
$\lim_{k \to \infty} \{ U(I^k) - e(A/I^k) \}  =  \infty$. This gives further evidence of the ``improved'' multiplicity  conjectures that suggest that \CM \ ideals with pure resolutions are the only ones for which the bounds are sharp.
\end{remark}

 \section{An Example}
In \cite{HEZH},  the authors state that its easy to construct examples of ideals with $\lim_{k\to
\infty} e(A/I^k) /U(I^k) = 1$. For sake of completeness we give a large class of ideals where 
$\lim_{k\to
\infty} e(A/I^k) /U(I^k) = 1$.  The notation will be as in section 3. Set $\m = (X_1,\ldots, X_n)$.

\s \label{exemp} Let $q \geq 2$ and let  $I \subseteq \m^q $ be a zero-dimensional ideal generated by $q$-forms.
It is easily verified that $\reg (I^k)  = qk + r$ for $k \gg 0$ (use \cite[3.2]{CHT}). Let 
$f_1,\ldots, f_n \in I$ be any regular sequence of $q$-forms. Set
$L = (f_1,\ldots, f_n)$.
Notice $e(L_\m, A_\m) = q^n = e(\m^q , A_\m)$. So by a theorem of Rees (see \ref{rees}), $L_\m$ is a reduction of $\m^qA_\m$.
It follows that $L$ is a reduction of $\m^q$.  Therefore  $L$ is also a reduction of 
$I$. By \ref{e-red} we get that $E(I) = E(L) = q^c$. So by \ref{formula} we get $\lim_{k\to
\infty} e(A/I^k) /U(I^k) = 1$.

\begin{remark}
We do not know as yet whether upper bound of multiplicity conjecture holds asymptotically for all ideals in the  class described
in \ref{exemp}. 
\end{remark}

\bibliographystyle{amsplain}

\end{document}